\documentclass[10pt]{amsart}
\usepackage{amsmath}
\usepackage{amsxtra}
\usepackage{amscd}
\usepackage{amsthm}
\usepackage{amsfonts}
\usepackage{amssymb}
\usepackage{eucal}
\usepackage[all]{xy}
\usepackage{graphicx}
\usepackage[usenames]{color}

\theoremstyle{remark}
\newtheorem{rem}[subsubsection]{Remark}

\theoremstyle{definition}

\theoremstyle{remark}

\newcommand{\secref}[1]{Sect.~\ref{#1}}

\numberwithin{equation}{section}

\newcommand{\nc}{\newcommand}
\nc{\renc}{\renewcommand}
\nc{\ssec}{\subsection}
\nc{\sssec}{\subsubsection}
\nc{\on}{\operatorname}

\nc{\ips}{{\iota_P^{(S)}}}
\nc{\ipms}{{\iota_{P^-}^{(S)}}}
\nc{\sfpps}{{\sfp_P^{(S)}}}
\nc{\sfppms}{{\sfp_{P^-}^{(S)}}}

\nc\ol{\overline}
\nc\wt{\widetilde}
\nc\tboxtimes{\wt{\boxtimes}}
\nc\tstar{\wt{\star}}
\nc{\alp}{\alpha}

\nc{\ZZ}{{\mathbb Z}}
\nc{\NN}{{\mathbb N}}
\nc{\OO}{{\mathbb O}}
\renc{\SS}{{\mathbb S}}
\nc{\DD}{{\mathbb D}}
\nc{\GG}{{\mathbb G}}

\nc{\Fq}{{\mathbb F}_q}
\nc{\Fqb}{\ol{\mathbb F}_q}
\nc{\Ql}{{\mathbb Q}_\ell}
\nc{\Qlb}{{\ol{\mathbb Q}_\ell}}
\nc{\id}{\text{id}}
\nc\X{\mathcal X}

\nc{\red}{\on{red}}
\nc{\Ho}{\on{Ho}}
\nc{\Hom}{\on{Hom}}
\nc{\coef}{\on{coef}}
\nc{\Lie}{\on{Lie}}
\nc{\Loc}{\on{Loc}}
\nc{\Pic}{\on{Pic}}
\nc{\Bun}{\on{Bun}}
\nc{\IC}{\on{IC}}
\nc{\Aut}{\on{Aut}}
\nc{\rk}{\on{rk}}
\nc{\Sh}{\on{Sh}}
\nc{\Perv}{\on{Perv}}
\nc{\pos}{{\on{pos}}}
\nc{\Conv}{\on{Conv}}
\nc{\Sph}{\on{Sph}}
\nc{\Sym}{\on{Sym}}
\nc{\BunBb}{\overline{\Bun}_B}
\nc{\BunNb}{\overline{\Bun}_N}
\nc{\BunTb}{\overline{\Bun}_T}
\nc{\BunBbm}{\overline{\Bun}_{B^-}}
\nc{\BunBbel}{\overline{\Bun}_{B,el}}
\nc{\BunBbmel}{\overline{\Bun}_{B^-,el}}
\nc{\Buno}{\overset{o}{\Bun}}
\nc{\BunPb}{{\overline{\Bun}_P}}
\nc{\BunBM}{\Bun_{B(M)}}
\nc{\BunBMb}{\overline{\Bun}_{B(M)}}
\nc{\BunPbw}{{\widetilde{\Bun}_P}}
\nc{\BunBP}{\widetilde{\Bun}_{B,P}}
\nc{\GUb}{\overline{G/U}}
\nc{\GUPb}{\overline{G/U(P)}}

\nc{\Hhom}{\underline{\on{Hom}}}
\nc\syminfty{\on{Sym}^{\infty}}
\nc\lal{\ol{\lambda}}
\nc\xl{\ol{x}}
\nc\thl{\ol{\theta}}
\nc\nul{\ol{\nu}}
\nc\mul{\ol{\mu}}
\nc\Sum\Sigma
\nc{\oX}{\overset{o}{X}{}}
\nc{\hl}{\overset{\leftarrow}h{}}
\nc{\hr}{\overset{\rightarrow}h{}}
\nc{\M}{{\mathcal M}}
\nc{\N}{{\mathcal N}}
\nc{\F}{{\mathcal F}}
\nc{\D}{{\mathcal D}}
\nc{\Q}{{\mathcal Q}}
\nc{\Y}{{\mathcal Y}}
\nc{\G}{{\mathcal G}}
\nc{\E}{{\mathcal E}}
\nc{\CalC}{{\mathcal C}}
\nc\Dh{\widehat{\D}}

\nc{\C}{{\mathcal C}}
\nc{\K}{{\mathcal K}}
\renewcommand{\H}{{\mathcal H}}

\nc{\T}{{\mathcal T}}
\nc{\V}{{\mathcal V}}
\renc{\P}{{\mathcal P}}
\nc{\A}{{\mathcal A}}
\nc{\B}{{\mathcal B}}
\nc{\U}{{\mathcal U}}

\nc{\Gr}{{\on{Gr}}}

\nc{\frn}{{\check{\mathfrak u}(P)}}

\nc{\fC}{\mathfrak C}
\nc{\p}{\mathfrak p}
\nc{\q}{\mathfrak q}
\nc\f{{\mathfrak f}}

\nc{\qo}{{\mathfrak q}}
\nc{\po}{{\mathfrak p}}
\nc{\s}{{\mathfrak s}}
\nc\w{\text{w}}

\nc\mathi\iota
\nc\Spec{\on{Spec}}
\nc\Proj{\on{Proj}}
\nc\Mod{\on{Mod}}
\nc{\tw}{\widetilde{\mathfrak t}}
\nc{\pw}{\widetilde{\mathfrak p}}
\nc{\qw}{\widetilde{\mathfrak q}}
\nc{\jw}{\widetilde j}

\nc{\grb}{\overline{\Gr}}
\nc{\I}{\mathcal I}

\nc{\lambdach}{{\check\lambda}}
\nc{\Lambdach}{{\check\Lambda}{}}
\nc{\much}{{\check\mu}}
\nc{\omegach}{{\check\omega}}
\nc{\nuch}{{\check\nu}}
\nc{\etach}{{\check\eta}}
\nc{\alphach}{{\check\alpha}}
\nc{\oblvtach}{{\check\oblvta}}
\nc{\rhoch}{{\check\rho}}
\nc{\ch}{{\check h}}

\nc{\Hb}{\overline{\H}}

\nc{\BA}{{\mathbb{A}}}
\nc{\BC}{{\mathbb{C}}}
\nc{\BG}{{\mathbb{G}}}
\nc{\BM}{{\mathbb{M}}}
\nc{\BO}{{\mathbb{O}}}
\nc{\BD}{{\mathbb{D}}}
\nc{\BN}{{\mathbb{N}}}
\nc{\BP}{{\mathbb{P}}}
\nc{\BQ}{{\mathbb{Q}}}
\nc{\BR}{{\mathbb{R}}}
\nc{\BZ}{{\mathbb{Z}}}
\nc{\BS}{{\mathbb{S}}}
\nc{\Deep}{{\bf{deep}}}
\nc{\deep}{deep}

\nc{\CA}{{\mathcal{A}}}
\nc{\CB}{{\mathcal{B}}}

\nc{\CE}{{\mathcal{E}}}
\nc{\CF}{{\mathcal{F}}}
\nc{\CH}{{\mathcal{H}}}

\nc{\CL}{{\mathcal{L}}}
\nc{\CC}{{\mathcal{C}}}
\nc{\CG}{{\mathcal{G}}}
\nc{\CalD}{{\mathcal{D}}}
\nc{\CM}{{\mathcal{M}}}
\nc{\CN}{{\mathcal{N}}}
\nc{\CK}{{\mathcal{K}}}
\nc{\CO}{{\mathcal{O}}}
\nc{\CP}{{\mathcal{P}}}
\nc{\CQ}{{\mathcal{Q}}}
\nc{\CR}{{\mathcal{R}}}
\nc{\CS}{{\mathcal{S}}}
\nc{\CT}{{\mathcal{T}}}
\nc{\CU}{{\mathcal{U}}}
\nc{\CV}{{\mathcal{V}}}
\nc{\CW}{{\mathcal{W}}}
\nc{\CX}{{\mathcal{X}}}
\nc{\CY}{{\mathcal{Y}}}
\nc{\CZ}{{\mathcal{Z}}}
\nc{\CI}{{\mathcal{I}}}

\nc{\csM}{{\check{\mathcal A}}{}}
\nc{\oM}{{\overset{\circ}{\mathcal M}}{}}
\nc{\obM}{{\overset{\circ}{\mathbf M}}{}}
\nc{\oCA}{{\overset{\circ}{\mathcal A}}{}}
\nc{\obA}{{\overset{\circ}{\mathbf A}}{}}
\nc{\ooM}{{\overset{\circ}{M}}{}}
\nc{\osM}{{\overset{\circ}{\mathsf M}}{}}
\nc{\vM}{{\overset{\bullet}{\mathcal M}}{}}
\nc{\nM}{{\underset{\bullet}{\mathcal M}}{}}
\nc{\oD}{{\overset{\circ}{\mathcal D}}{}}
\nc{\obD}{{\overset{\circ}{\mathbf D}}{}}
\nc{\oA}{{\overset{\circ}{\mathbb A}}{}}
\nc{\op}{{\overset{\bullet}{\mathbf p}}{}}
\nc{\cp}{{\overset{\circ}{\mathbf p}}{}}
\nc{\oU}{{\overset{\bullet}{\mathcal U}}{}}
\nc{\oZ}{{\overset{\circ}{\mathcal Z}}{}}
\nc{\ofZ}{{\overset{\circ}{\mathfrak Z}}{}}
\nc{\oF}{{\overset{\circ}{\fF}}}

\nc{\fa}{{\mathfrak{a}}}
\nc{\fb}{{\mathfrak{b}}}
\nc{\fd}{{\mathfrak{d}}}
\nc{\ff}{{\mathfrak{f}}}
\nc{\fg}{{\mathfrak{g}}}
\nc{\fgl}{{\mathfrak{gl}}}
\nc{\fh}{{\mathfrak{h}}}
\nc{\fj}{{\mathfrak{j}}}
\nc{\fl}{{\mathfrak{l}}}
\nc{\fm}{{\mathfrak{m}}}
\nc{\fn}{{\mathfrak{n}}}
\nc{\fu}{{\mathfrak{u}}}
\nc{\fp}{{\mathfrak{p}}}
\nc{\fr}{{\mathfrak{r}}}
\nc{\fs}{{\mathfrak{s}}}
\nc{\ft}{{\mathfrak{t}}}
\nc{\fz}{{\mathfrak{z}}}
\nc{\fsl}{{\mathfrak{sl}}}
\nc{\hsl}{{\widehat{\mathfrak{sl}}}}
\nc{\hgl}{{\widehat{\mathfrak{gl}}}}
\nc{\hg}{{\widehat{\mathfrak{g}}}}
\nc{\chg}{{\widehat{\mathfrak{g}}}{}^\vee}
\nc{\hn}{{\widehat{\mathfrak{n}}}}
\nc{\chn}{{\widehat{\mathfrak{n}}}{}^\vee}

\nc{\fA}{{\mathfrak{A}}}
\nc{\fB}{{\mathfrak{B}}}
\nc{\fD}{{\mathfrak{D}}}
\nc{\fE}{{\mathfrak{E}}}
\nc{\fF}{{\mathfrak{F}}}
\nc{\fG}{{\mathfrak{G}}}
\nc{\fK}{{\mathfrak{K}}}
\nc{\fL}{{\mathfrak{L}}}
\nc{\fM}{{\mathfrak{M}}}
\nc{\fN}{{\mathfrak{N}}}
\nc{\fP}{{\mathfrak{P}}}
\nc{\fU}{{\mathfrak{U}}}
\nc{\fV}{{\mathfrak{V}}}
\nc{\fZ}{{\mathfrak{Z}}}

\nc{\bb}{{\mathbf{b}}}
\nc{\bc}{{\mathbf{c}}}
\nc{\bd}{{\mathbf{d}}}
\nc{\bbf}{{\mathbf{f}}}
\nc{\be}{{\mathbf{e}}}
\nc{\bi}{{\mathbf{i}}}
\nc{\bj}{{\mathbf{j}}}
\nc{\bn}{{\mathbf{n}}}
\nc{\bo}{{\mathbf{o}}}
\nc{\bp}{{\mathbf{p}}}
\nc{\bq}{{\mathbf{q}}}
\nc{\bu}{{\mathbf{u}}}
\nc{\bv}{{\mathbf{v}}}
\nc{\bx}{{\mathbf{x}}}
\nc{\bs}{{\mathbf{s}}}
\nc{\by}{{\mathbf{y}}}
\nc{\bw}{{\mathbf{w}}}
\nc{\bA}{{\mathbf{A}}}
\nc{\bK}{{\mathbf{K}}}
\nc{\bB}{{\mathbf{B}}}
\nc{\bC}{{\mathbf{C}}}
\nc{\bG}{{\mathbf{G}}}
\nc{\bD}{{\mathbf{D}}}
\nc{\bH}{{\mathbf{H}}}
\nc{\bM}{{\mathbf{M}}}
\nc{\bN}{{\mathbf{N}}}
\nc{\bO}{{\mathbf{O}}}
\nc{\bV}{{\mathbf{V}}}
\nc{\bW}{{\mathbf{W}}}
\nc{\bX}{{\mathbf{X}}}
\nc{\bZ}{{\mathbf{Z}}}
\nc{\bS}{{\mathbf{S}}}

\nc{\sA}{{\mathsf{A}}}
\nc{\sB}{{\mathsf{B}}}
\nc{\sC}{{\mathsf{C}}}
\nc{\sD}{{\mathsf{D}}}
\nc{\sF}{{\mathsf{F}}}
\nc{\sG}{{\mathsf{G}}}
\nc{\sK}{{\mathsf{K}}}
\nc{\sM}{{\mathsf{M}}}
\nc{\sO}{{\mathsf{O}}}
\nc{\sW}{{\mathsf{W}}}
\nc{\sQ}{{\mathsf{Q}}}
\nc{\sP}{{\mathsf{P}}}
\nc{\sR}{{\mathsf{R}}}
\nc{\sZ}{{\mathsf{Z}}}
\nc{\sfp}{{\mathsf{p}}}
\nc{\sfq}{{\mathsf{q}}}
\nc{\sr}{{\mathsf{r}}}
\nc{\bk}{{\mathsf{k}}}
\nc{\sg}{{\mathsf{g}}}
\nc{\sff}{{\mathsf{f}}}
\nc{\sfb}{{\mathsf{b}}}
\nc{\sfc}{{\mathsf{c}}}
\nc{\sd}{{\mathsf{d}}}

\nc{\BK}{{\bar{K}}}

\nc{\tA}{{\widetilde{\mathbf{A}}}}
\nc{\tB}{{\widetilde{\mathcal{B}}}}
\nc{\tg}{{\widetilde{\mathfrak{g}}}}
\nc{\tG}{{\widetilde{G}}}
\nc{\TM}{{\widetilde{\mathbb{M}}}{}}
\nc{\tO}{{\widetilde{\mathsf{O}}}{}}
\nc{\tU}{{\widetilde{\mathfrak{U}}}{}}
\nc{\TZ}{{\tilde{Z}}}
\nc{\tx}{{\tilde{x}}}
\nc{\tbv}{{\tilde{\bv}}}
\nc{\tfP}{{\widetilde{\mathfrak{P}}}{}}
\nc{\tz}{{\tilde{\zeta}}}
\nc{\tmu}{{\tilde{\mu}}}

\nc{\urho}{\underline{\rho}}
\nc{\uB}{\underline{B}}
\nc{\uC}{{\underline{\mathbb{C}}}}
\nc{\ui}{\underline{i}}
\nc{\uj}{\underline{j}}
\nc{\ofP}{{\overline{\mathfrak{P}}}}
\nc{\oB}{{\overline{\mathcal{B}}}}
\nc{\og}{{\overline{\mathfrak{g}}}}
\nc{\oI}{{\overline{I}}}

\nc{\eps}{\varepsilon}
\nc{\hrho}{{\hat{\rho}}}

\nc{\one}{{\mathbf{1}}}
\nc{\two}{{\mathbf{t}}}

\nc{\Rep}{{\mathop{\operatorname{\rm Rep}}}}
\nc{\Tot}{{\mathop{\operatorname{\rm Tot}}}}
\nc{\Ker}{{\mathop{\operatorname{\rm Ker}}}}
\nc{\im}{{\mathop{\operatorname{\rm Im}}}}
\nc{\Hilb}{{\mathop{\operatorname{\rm Hilb}}}}
\nc{\End}{{\mathop{\operatorname{\rm End}}}}
\nc{\Ext}{{\mathop{\operatorname{\rm Ext}}}}
\nc{\CHom}{{\mathop{\operatorname{{\mathcal{H}}\it om}}}}
\nc{\GL}{{\mathop{\operatorname{\rm GL}}}}
\nc{\gr}{{\mathop{\operatorname{\rm gr}}}}
\nc{\HN}{{\mathop{\operatorname{\rm HN}}}}
\nc{\Id}{{\mathop{\operatorname{\rm Id}}}}
\nc{\de}{{\mathop{\operatorname{\rm def}}}}
\nc{\length}{{\mathop{\operatorname{\rm length}}}}
\nc{\supp}{{\mathop{\operatorname{\rm supp}}}}

\nc{\Cliff}{{\mathsf{Cliff}}}
\nc{\Fl}{\on{Fl}}
\nc{\Fib}{{\mathsf{Fib}}}
\nc{\Coh}{{\on{Coh}}}
\nc{\QCoh}{{\on{QCoh}}}
\nc{\IndCoh}{{\on{IndCoh}}}
\nc{\FCoh}{{\mathsf{FCoh}}}

\nc{\reg}{{\text{\rm reg}}}

\nc{\cplus}{{\mathbf{C}_+}}
\nc{\cminus}{{\mathbf{C}_-}}
\nc{\cthree}{{\mathbf{C}_\bullet}}
\nc{\Qbar}{{\bar{Q}}}
\nc\Eis{\on{Eis}}
\nc\Eisb{\ol\Eis{}}
\nc\Eisr{\on{Eis}^{rat}{}}
\nc\wh{\widehat}
\nc{\Def}{\on{Def_{\check{\fb}}(E)}}
\nc{\barZ}{\overline{Z}{}}
\nc{\barbarZ}{\overline{\barZ}{}}
\nc{\barpi}{\overline\pi}
\nc{\barbarpi}{\overline\barpi}
\nc{\barpip}{\overline\pi{}^+}
\nc{\barpim}{\overline\pi{}^-}

\nc{\fq}{\mathfrak q}

\nc{\fqb}{\ol{\sfq}{}}
\nc{\fpb}{\ol{\sfp}{}}
\nc{\fpr}{{\sfp^{rat}}{}}
\nc{\fqr}{{\sfq^{rat}}{}}

\nc{\hattimes}{\wh\otimes}

\nc{\bh}{{\bar{h}}}
\nc{\bOmega}{{\overline{\Omega(\check \fn)}}}

\nc{\seq}[1]{\stackrel{#1}{\sim}}

\nc{\cT}{{\check{T}}}
\nc{\cG}{{\check{G}}}
\nc{\cM}{{\check{M}}}
\nc{\cB}{{\check{B}}}

\nc{\ct}{{\check{\mathfrak t}}}
\nc{\cg}{{\check{\fg}}}
\nc{\cb}{{\check{\fb}}}
\nc{\cn}{{\check{\fn}}}

\nc{\cLambda}{{\check\Lambda}}

\nc{\cla}{{\check\lambda}}
\nc{\cmu}{{\check\mu}}
\nc{\cnu}{{\check\nu}}
\nc{\ceta}{{\check\eta}}

\nc{\DefbE}{{\on{Def}_{\cB}(E_\cT)}}

\nc{\imathb}{{\ol{\imath}}}
\nc{\rlr}{\overset{\longrightarrow}{\underset{\longrightarrow}\longleftarrow}}

\nc{\oBun}{\overset{\circ}\Bun}
\nc{\LocSys}{\on{LocSys}}
\nc{\BunBbb}{\ol{\ol{Bun}}_B}
\nc{\BunBr}{\Bun_B^{rat}}
\nc{\BunBrsg}{\Bun_B^{rat,\on{s.g.}}}
\nc{\BunBrp}{\Bun_B^{rat,polar}}
\nc{\BunBrpbg}{\Bun_B^{rat,polar,\on{b.g.}}}
\nc{\BunBrpsg}{\Bun_B^{rat,polar,\on{s.g.}}}
\nc{\BunTrp}{\Bun_T^{rat,polar}}
\nc{\BunTrpbg}{\Bun_T^{rat,polar,\on{b.g.}}}
\nc{\BunTrpsg}{\Bun_T^{rat,polar,\on{s.g.}}}
\nc{\BunNr}{\Bun_N^{rat}}
\nc{\BunNre}{\Bun_N^{enh,rat}}
\nc{\BunTr}{\Bun_T^{rat}}
\nc{\Vect}{\on{Vect}}
\nc{\Whit}{\on{Whit}}
\nc{\CTb}{\ol{\on{CT}}}
\nc{\Ran}{\on{Ran}}
\nc{\CTr}{\on{CT}^{rat}{}}
\nc\jmathr{\jmath^{rat}{}}
\nc{\ux}{\underline{x}}
\nc{\clambda}{{\check\lambda}}
\nc{\calpha}{{\check\alpha}}
\nc{\ind}{{\mathbf{ind}}}
\nc{\oblv}{{\mathbf{oblv}}}
\nc{\ox}{{\overline{x}}}
\nc{\cLa}{\check{\Lambda}}
\nc{\StinftyCat}{\on{DGCat}}
\nc{\inftyCat}{\infty\on{-Cat}}
\nc{\inftygroup}{\infty\on{-Grpd}}
\nc{\Dmod}{\on{D-mod}}
\nc{\CMaps}{{\mathcal Maps}}
\nc{\Maps}{\on{Maps}}
\nc{\affSch}{\on{Sch}^{\on{aff}}}
\nc{\dr}{{\on{dR}}}
\nc{\oCY}{\overset{\circ}\CY}
\nc{\leqG}{\underset{G}\leq}
\nc{\leqM}{\underset{M}\leq}
\nc{\leqGad}{\underset{G_{ad}}\leq}
\nc{\leqMad}{\underset{M_{ad}}\leq}
\nc{\Tr}{\on{Tr}}
\nc{\Frob}{\on{Frob}}
\nc{\DGCat}{\on{DGCat}}
\nc{\tDGCat}{2\on{-DGCat}}
\nc{\ev}{\on{ev}}
\nc{\mmod}{\on{-}\mathbf{mod}}
\nc{\sotimes}{\overset{!}\otimes}

\begin{document}

\title[From geometric to function-theoretic Langlands]{From geometric to function-theoretic Langlands \\  (or how to invent shtukas)}
\author{Dennis Gaitsgory}

\begin{abstract}
This is an informal note that explains that the classical Langlands theory over function fields can be obtained from the geometric
one by taking the trace of Frobenius. The operation of taking the trace of Frobenius takes place at the categorical level,
and this we deduce that the \emph{space} of automorphic functions is the trace of the Frobenius on the \emph{category}
of automorphic sheaves.
\end{abstract}

\date{\today}

\maketitle

\section{Trace of Frobenius on a category}

\ssec{The trace of an endomorphism}

\sssec{}  \label{sss:trace}

The notion of a trace of an endomorphism of a dualizable object $\bo$ in a symmetric monoidal
category $\bO$ is well-known:

\medskip

If $\bo\in \bO$ is dualizable and $T:\bO\to \bO$ consider the corresponding map $$Q_T:\one_\bO\to \bo\otimes \bo^\vee$$
and $\Tr(T,\bo)\in \End_\bO(\one_\bO)$ is by definition the composition 
$$\one_\bO\overset{Q_T}\to \bo\otimes \bo^\vee\overset{\on{ev}}\longrightarrow \one_\bO.$$

\sssec{}  \label{sss:endo of stack}

Set $\bO=\DGCat$, and let $\bo=\Dmod(\CY)$, where is a quasi-compact algebraic stack with affine diagonal. It is known that
$\Dmod(\CY)$ is dualizable (as a DG category), where the evaluation functor
$$\Dmod(\CY)\otimes \Dmod(\CY)\to \Vect$$
is given by
$$\Dmod(\CY)\otimes \Dmod(\CY)\overset{\boxtimes}\simeq \Dmod(\CY\times \CY)\overset{\Delta_\CY^!} \longrightarrow
\Dmod(\CY)\overset{\on{H}^\cdot(\CY,\omega_\CY)}\to \Vect,$$
where $\on{H}^\cdot(\CY,-)$ is functor of sheaf cohomology (refined to the chain level). 

\medskip

Suppose that $T:\Dmod(\CY)\to \Dmod(\CY)$
is given by pullback with respect to an endomorphism $F:\CY\to \CY$. 

\medskip

Then the claim is that 
$$\Tr(F^!,\Dmod(\CY))\simeq \on{H}^\cdot(\CY^F,\omega_{\CY^F})=\on{H}_\cdot^{\on{BM}}(\CY^F),$$
where $\CY^F$ is the stack-theoretic fixed-point locus of $F$, i.e., 
$$\CY\underset{\Gamma_F,\CY\times \CY,\Delta_\CY}\times \CY,$$
and where $\on{H}_\cdot^{\on{BM}}(-)$ stands for Borel-Moore chains.

\medskip

Indeed, this follows by base change from the fact that the corresponding to $Q_T$ is given by
$$\Vect \overset{\omega_\CY}\longrightarrow \Dmod(\CY) \overset{(\on{Graph}_F)_\bullet}\longrightarrow \Dmod(\CY\times \CY) \simeq 
\Dmod(\CY)\otimes \Dmod(\CY).$$

\sssec{}   \label{sss:endo of constr}

Let us now change the context, where $\CY$ is a stack over $\Fqb$. We let $D(\CY)$ be the ind-completion usual bounded 
constructible derived category. 

\medskip

The problem is that in the constructible setting, the functor 
$$D(\CY_1)\otimes D(\CY_2)\to D(\CY_1\times \CY_2)$$
is \emph{no longer} an equivalence. However, we will pretend that it is. That said, Drinfeld had an idea how to provide
a framework for this; I think this amounts to tweaking the definition of $\DGCat$. 

\medskip

When considering functors $D(\CY_1)\to D(\CY_2)$, we will restrict ourselves to functors given by kernels: i.e., to
an object $Q\in D(\CY_1\times \CY_2)$ we attach a functor
$$\CF\mapsto (\on{pr}_{\CY_2})_\bullet((\on{pr}_{\CY_1})^!\otimes Q.$$

When talking about a functor admitting a right adjoint, we will also mean that this functor is given by a kernel. 

\sssec{}  \label{sss:Frob}

Assume now that $\CY$ is defined over $\Fq$, and $F$ is the corresponding geometric Frobenius
endomorphism $\Frob_\CY$ of $\CY$.  Then the calculation of \secref{sss:endo of stack} implies that
\begin{equation} \label{e:space of functions from sheaves}
\Tr(\Frob^!_\CY,D(\CY))\simeq \on{Funct}(\CY(\Fq),\Qlb).
\end{equation}

\ssec{Functoriality}

\sssec{}  \label{sss:trace of 1-morph}

Suppose that in the context of \secref{sss:trace}, $\bO$ is actually a monoidal 2-category. Let
$$(\bo,T_\bo:\bo\to \bo) \text{ and } (\bo',T_{\bo'}:\bo\to \bo')$$ be two pairs of objects, each 
equipped with an endomorphism.
Let
$$S:\bo\to \bo'$$ be a 1-morphism, equipped with a 2-morphism
\begin{equation} \label{e:2-morph}
S\circ T_\bo\to T_{\bo'}\circ S.
\end{equation} 

Assume also that $S$ \emph{admits a right adjoint} in $\bO$ (this is an intrinsic 2-categorical condition). 

\medskip

We claim that in this case there is a canonical map
$$\Tr(S):\Tr(T_\bo,\bo)\to \Tr(T_{\bo'},{\bo'})$$
in the 1-category $\End(\one_\bO)$. 

\medskip

Indeed,  $\Tr(S)$ is given as a composition
\begin{multline*}
\ev_\bo\circ Q_{T_\bo}\to \ev_{\bo'}\circ (S\otimes (S^R)^\vee) \circ Q_{T_\bo} \simeq 
\ev_{\bo'}\circ (\on{id}_{\bo'}\otimes (S^R)^\vee) \circ Q_{S\circ T_\bo} \to \\
\to \ev_{\bo'}\circ (\on{id}_{\bo'}\otimes (S^R)^\vee) \circ Q_{T_{\bo'}\circ S}\simeq 
\ev_{\bo'}\circ Q_{T_{\bo'}\circ S\circ S^R}\to \ev_{\bo'}\circ Q_{T_{\bo'}}.
\end{multline*} 

\sssec{}  \label{sss:comp}

One checks directly that the above construction is compatible with compositions. I.e., if
we have three pairs $(\bo,T_\bo:\bo\to \bo)$, $(\bo,T_{\bo'}:\bo'\to \bo')$ and $(\bo'',T_{\bo''}:\bo''\to \bo'')$
and the 1-mophisms 
$$S':\bo\to \bo', \quad S'':\bo'\to \bo'',$$
and also the corresponding 2-morphisms \eqref{e:2-morph}, then
$$\Tr(S''\circ S')\simeq \Tr(S'')\circ \Tr(S'),$$
as morphisms in $\End(\one_\bO)$. 

\sssec{}

For example, take $\bO=\DGCat$ and let $(\bo,T)$ be an object with an endomorphism, i.e., a DG category $\CC$ with
an endo-functor $T$. Let $c\in \CC$ be a \emph{compact} object, equipped with a morphism
$$c\to T(c).$$

We can view such $c$ as a datum of 1-morphism
$$\Vect=\one_{\DGCat}\to \CC$$
that admits a (continuous!) right adjoint and a 2-morphism as in \eqref{e:2-morph}. 
 
\medskip

The construction in \secref{sss:trace of 1-morph} yields a map in $\Vect$
$$k\to \Tr(T,\CC),$$
where $k$ is the ground field (i.e., the unit object in $\Vect$). I.e., we obtain an element in the vector
space $\Tr(T,\CC)$; we denote this element by $\Tr(T,c)$. 

\sssec{}

Let us return to the example of \secref{sss:Frob}. Let $\CM$ be a compact object in $D(\CY)$, equipped with a map
\begin{equation} \label{e:Weil}
\CM\to \Frob_\CY^!(\CM).
\end{equation}

One checks that
$$\Tr(\Frob_\CY^!,\CM)\in \Tr(\Frob^!_\CY,D(\CY))\simeq \on{Funct}(\CY(\Fq),\Qlb)$$
is the same as the function obtained from $\CM$ by the usual faisceaux-fonctions. 

\begin{rem}
Here we are using the following version of faisceaux-fonctions: to $\CM\in D(\CY)$, equipped with a map \eqref{e:Weil},
we associate the function of $\CY(\Fq)$ equal to traces of the Frobenius on the !-fibers of $\CM$. This equals the usual
faisceaux-fonctions (i.e., for *-fibers) on the Verdier dual of $\CM$.
\end{rem} 

\ssec{Action of local systems}  \label{ss:act loc sys}

\sssec{}  \label{sss:trace QCoh}

In the general context of \secref{sss:trace} let us again take $\bO=\DGCat$ but $\bo=\QCoh(\CY)$. Note that
$\QCoh(\CY)$ is again self-dual with the evaluation functor being
$$\QCoh(\CY)\otimes \QCoh(\CY)\simeq \QCoh(\CY\times \CY) \overset{\Delta^*}\longrightarrow 
\QCoh(\CY) \overset{\Gamma(\CY,-)}\longrightarrow \Vect.$$

\medskip

Let $T$ be again given by pullback along an endomorphism $F$. Then the same calculation as in \secref{sss:endo of stack}
shows that
$$\Tr(F^*,\QCoh(\CY))\simeq \Gamma(\CY^T,\CO_{\CY^T}).$$

\medskip

By the functoriality developed in \secref{sss:trace of 1-morph}, the structure of symmetric monoidal category on 
$\QCoh(\CY)$ defines a structure of commutative algebra on $\Gamma(\CY^T,\CO_{\CY^T})$. It is straightforward to
check that this is the usual structure of commutative algebra on $\Gamma(\CY^T,\CO_{\CY^T})$. 

\sssec{}  \label{sss:Vincent action}

We want to take $\CY$ to be the stack $\LocSys$ of \'etale local systems on a curve over $\Fqb$. Of course, it does not exist as an algebraic
stack. But we will pretend that it does. There is a hope that the corresponding monoidal category $\QCoh(\LocSys)$, or whatever we need from
it, could actually be defined, as an algebra object in $\DGCat$. 

\medskip

If our curve is defined over $\Fq$, then $\LocSys$ acquires an automorphism, given by Frobenius. By \secref{sss:trace QCoh},
the corresponding category
$$\Tr(\Frob_{\LocSys}^*,\QCoh(\LocSys))$$
identifies with the vector space $\Gamma(\LocSys^{\on{arthm}},\CO_{\LocSys^{\on{arthm}}})$, where 
$\LocSys^{\on{arthm}}$ is the stack of \emph{arithmetic} local systems. 

\sssec{}

We will now assume the \emph{geometric spectral decomposition}, i.e., the action of the monoidal category
$\QCoh(\LocSys)$ on $D(\Bun_G)$.
Recall that such an action does indeed exist in the context of D-modules, by the ``generalized vanishing theorem",
see \cite[Corollary 4.5.5]{Ga}.

\medskip

Applying the functoriality construction from \secref{sss:trace of 1-morph}, we obtain an action of the algebra 
$\Gamma(\LocSys^{\on{arthm}},\CO_{\LocSys^{\on{arthm}}})$ on $\on{Funct}(\Bun_G(\Fq),\Qlb)$.

\medskip

We claim that this is the action constructed in Vincent Lafforgue's work \cite{Laf}, using shtukas. A convenient way to do this is
to first reinterpret Lafforgue's construction \`a la Drinfeld, i.e., organize the cohomologies of shtukas into
an object of $\QCoh(\LocSys^{\on{arthm}})$. This will be done in the next section.

\ssec{Hecke action}

\sssec{}

Fix a rational point $x\in X$. For a representation $V$ of the dual group $\cG$, we have the Hecke functor
$$\on{H}_{x,V}:D(\Bun_G)\to D(\Bun_G),$$
which is naturally compatible with the Frobenius endo-functor $\Frob^!_{\Bun_G}$ on both sides.

\medskip

By \secref{sss:trace of 1-morph}, it gives rise to an endomorphism
$$\on{Tr}(\on{H}_{x,V}): \on{Funct}(\Bun_G(\Fq),\Qlb)\to \on{Funct}(\Bun_G(\Fq),\Qlb).$$

It should be possible to see (but I haven't done that yet) that this endomorphism is the usual Hecke
functor corresponding to $V$ at $x$. 

\sssec{}

Restriction to the formal disc around $x$ defines a map
$$\LocSys^{\on{arthm}}\to \cG/\on{Ad}(\cG),$$
where we think of $\cG/\on{Ad}(\cG)$ as the stack of unramified arithmetic local systems on the disc around $x$.

\medskip

In particular, we obtain a map
$$\CH^{\on{cl}}\simeq \Gamma(\cG/\on{Ad}(\cG),\CO_{\cG/\on{Ad}(\cG)})\to \Gamma(\LocSys^{\on{arthm}},\CO_{\LocSys^{\on{arthm}}}),$$
where $\CH^{\on{cl}}$ is the classical spherical Hecke algebra.

\medskip

Combining with \secref{sss:Vincent action} we obtain an action of $\CH^{\on{cl}}$ on $\on{Funct}(\Bun_G(\Fq),\Qlb)$. 

\medskip

It should be possible to see (but I haven't done that either) that for a representation $V$ of the dual group, the action 
of the corresponding element $\on{H}^{\on{cl}}_{x,V}\in \CH^{\on{cl}}$ on $\on{Funct}(\Bun_G(\Fq),\Qlb)$ equals 
$\on{Tr}(\on{H}_{x,V})$. 

\medskip

If this is the case, this gives a ``conceptual" explanation of Vincent's main formula that expresses the Hecke operators 
as particular excursion operators. 

\section{Enhanced trace}

\ssec{Trace on 2-categories}

\sssec{}

We now want to take our monoidal 2-category $\bO$ to be that of \emph{2-categories tensored over $\DGCat$},
denoted $\tDGCat$. We do not quite know how to give a completely satisfactory definitions, but modulo the questions 
of 1-affineness, the following will do:

\medskip

We let the objects of $\tDGCat$ be monoidal DG categories (i.e., associative algebra objects in $\DGCat$). For two
such, denoted $\CA_0$ and $\CA_1$, we let the category of 1-morphisms $\CA_0\to \CA_1$ to be that of 
$(\CA_1,\CA_0)$-bimodules. Note that the latter is naturally a 2-category, but for our purposes we will not
need to consider non-invertible 2-morphisms in it. 

\medskip

In other words, we are thinking of a monoidal DG category $\CA$ in terms of the 2-category $\CA\mmod$. 
And the category of $(\CA_1,\CA_0)$-bimodules is naturally that of functors
$$\CA_0\mmod\to \CA_1\mmod.$$

\medskip

So, we will denote objects of $\tDGCat$ by $\CA\mmod$ (rather than $\CA$). 

\sssec{}

The symmetric monoidal structure on $\tDGCat$ is given by 
$$\CA,\CB\mapsto \CA\otimes \CB,$$
i.e.,
$$\CA\mmod\otimes \CB\mmod:=(\CA\otimes \CB)\mmod.$$

\medskip

The unit object in $\tDGCat$ is $\CA=\Vect$ so that
$$\CA\mmod=\DGCat.$$

\medskip

Note that
$$\End_{\tDGCat}(\one_{\tDGCat})=\DGCat$$
as a symmetric monoidal category. 

\sssec{}  \label{sss:higher trace}

Let $\CA\mmod$ be an object of $\tDGCat$, and let $\CT$ be its endofunctor, i.e., an $\CA$-bimodule category. Then
$$\Tr(\CT,\CA\mmod)=\CT\underset{\CA^{\on{rev-mult}}\otimes \CA}\otimes \CA\in \DGCat,$$
I.e., this Hochschild homology of $\CA$ with coefficients in $\CT$, which we will also denote by $\on{HH}_\cdot(\CA,\CT)$.

\medskip

Let $\CC$ be an object of $\CA\mmod$. Assume that $\CC$ is \emph{dualizable as an $\CA$-module category}. This is equivalent to the
condition that the 1-morphism
$$\DGCat\to \CA\mmod$$ 
in $\tDGCat$ admits a right adjoint. 

\medskip

Assume that we are given a functor of $\CA$-module categories 
$$T:\CC\to \CT\underset{\CA}\otimes \CC.$$

The the construction of \secref{sss:trace of 1-morph} produces an object 
$$\Tr_{\CT,\CA}(T,\CC)\in \Tr(\CT,\CA\mmod).$$

\ssec{Drinfeld's $\CO$-module}  \label{ss:Drinf}

\sssec{}

Let us take $\CA=\QCoh(\LocSys)$, with $\CT$ given by the pullback with respect to $\Frob_{\LocSys}$. Note that the DG category
$$\Tr(\Frob_{\LocSys}^*,\QCoh(\LocSys))$$
is by definition 
$$\QCoh(\LocSys)\underset{\on{Graph}_{\Frob_{\LocSys}},\QCoh(\LocSys)\otimes \QCoh(\LocSys),\Delta_{\LocSys}}\otimes \QCoh(\LocSys),$$
while the latter identifies with 
$$\QCoh(\LocSys^{\on{arthm}}),$$
i.e.,
$$\Tr(\Frob_{\LocSys}^*,\QCoh(\LocSys))\simeq \QCoh(\LocSys^{\on{arthm}}).$$

\sssec{}

Let us now take $\CC=D(\Bun_G)$, regarded as an $\QCoh(\LocSys)$-module category. We let the datum of $T$ be given by $\Frob^!_{\Bun_G}$.
Applying the construction of \secref{sss:higher trace}, we obtain an object
$$\Tr_{\Frob_{\LocSys}^*,\QCoh(\LocSys)}(\Frob^!_{\Bun_G},D(\Bun_G))=:\on{Drinf}\in \QCoh(\LocSys^{\on{arthm}}).$$

\sssec{}

Note that the compatibility with compositions in \secref{sss:comp} implies that
$$\Gamma(\LocSys^{\on{arthm}},\on{Drinf})
\simeq \Tr(\Frob^!_{\Bun_G},D(\Bun_G)),$$
while the latter identifies 
$$\on{Funct}(\Bun_G(\Fq),\Qlb).$$

\sssec{}

Let us now take $\CC=D(\Bun_G)_{\on{temp}}$, and denote
$$\Tr_{\Frob_{\LocSys}^*,\QCoh(\LocSys)}(\Frob^!_{\Bun_G},D(\Bun_G)_{\on{temp}})=:\on{Drinf}_{\on{temp}}.$$

If we believe that the action of $\QCoh(\LocSys)$ on the Whittaker object defines an equivalence
$$\QCoh(\LocSys)\to D(\Bun_G)_{\on{temp}},$$
we obtain an equivalence
$$\on{Drinf}_{\on{temp}}\simeq \CO_{\LocSys^{\on{arthm}}}.$$

\ssec{Relation to shtukas}  \label{ss:sht}

\sssec{}

We will now show that the object $\on{Drinf}\in \QCoh(\LocSys^{\on{arthm}})$ is isomorphic to the 
one arising from cohomologies of shtukas, denoted $\on{Drinf-Sht}$.  

\medskip

The object $\on{Drinf-Sht}$ was characterized by the following property. Fix a rational point $x\in X$. Let
$\CE_{V,x}$ be the vector bundle on $\LocSys$ associated to $x$ and a representation $V$ of $\cG$. 

\medskip

Then
$$\Gamma(\LocSys^{\on{arthm}},\CE_{V,x}|_{\LocSys^{\on{arthm}}}\otimes \on{Drinf-Sht})$$
is by definition
$$\on{H}^\cdot(\Bun_G\underset{\Frob_{\Bun_G},\Bun_G\times \Bun_G}\times \CH_x,\CS_V),$$
where $\CH_x$ is the Hecke stack at $x$ and $\CS_V\in D(\CH_x)$ is the object corresponding to $V$ by
geometric Satake.

\medskip

In order to establish the desired isomorphism $\on{Drinf}\simeq \on{Drinf-Sht}$, we will 
construct an isomorphism 
\begin{equation} \label{e:key}
\Gamma(\QCoh(\LocSys^{\on{arthm}},\CE_{V,x}|_{\LocSys^{\on{arthm}}}\otimes \on{Drinf})\simeq 
\on{H}^\cdot(\Bun_G\underset{\Frob_{\Bun_G},\Bun_G\times \Bun_G}\times \CH_x,\CS_V)
\end{equation} 

\begin{rem}

By making $x$ move along $X$, the isomorphism
$$\Gamma(\QCoh(\LocSys^{\on{arthm}},\CE_{V,x}|_{\LocSys^{\on{arthm}}}\otimes \on{Drinf-Sht})\simeq 
\on{H}(\Bun_G\underset{\Frob_{\Bun_G},\Bun_G\times \Bun_G}\times \CH_x,\CS_V)$$
extends to one between Weil sheaves on $X$.

\medskip

Hopefully, the same will be the case for the isomorphism \eqref{e:key} that we are about to construct.

\end{rem} 

\sssec{}

In order to establish \eqref{e:key} let us unravel the definition of the left-hand side. It is the composition
\begin{multline}  \label{e:long composition}
\Vect \to D(\Bun_G)\underset{\QCoh(\LocSys)}\otimes D(\Bun_G) \simeq \\
\simeq \left(D(\Bun_G)\otimes D(\Bun_G)\right)\underset{\QCoh(\LocSys\times \LocSys)}\otimes \QCoh(\LocSys)\to \\
\QCoh(\LocSys)\underset{\on{Graph}_{\Frob_{\LocSys}},\QCoh(\LocSys\times \LocSys),\Delta_{\LocSys}}\otimes \QCoh(\LocSys)\simeq \\
\simeq \QCoh(\LocSys^{\on{arthm}}) 
\overset{\otimes \CE_{V,x}|_{\LocSys^{\on{arthm}}}}\longrightarrow \QCoh(\LocSys^{\on{arthm}}) 
\overset{\Gamma(\LocSys^{\on{arthm}},-)}\longrightarrow \Vect, 
\end{multline}
where the first arrow 
$$\Vect \to D(\Bun_G)\underset{\QCoh(\LocSys)}\otimes D(\Bun_G),$$
corresponds to the unit of the self-duality datum for $D(\Bun_G)$ as a module category over $\QCoh(\LocSys)$, 
and the third arrow is induced by the arrow 
$$\Phi:D(\Bun_G)\otimes D(\Bun_G)\to \QCoh(\LocSys),$$
defined as follows: for $\CE\in \QCoh(\LocSys)$ and $\CM_1,\CM_2\in \Dmod(\Bun_G)$, we have
\begin{equation}   \label{e:interp counit} 
\Gamma(\LocSys,\CE\otimes \Phi(\CM_1,\CM_2))=\on{H}^\cdot\left(\Bun_G,\Frob^!_{\Bun_G}(\CM_1) \sotimes (\CE\star \CM_2)\right).
\end{equation}

\sssec{}

Note that the last three lines in \eqref{e:long composition} can be replaced by 
\begin{multline*} 
\left(D(\Bun_G)\otimes D(\Bun_G)\right)\underset{\QCoh(\LocSys\times \LocSys)}\otimes \QCoh(\LocSys)\to  \\
\to D(\Bun_G)\otimes D(\Bun_G) \overset{\Phi}\longrightarrow 
\QCoh(\LocSys)\overset{\otimes \CE_{V,x}}\longrightarrow \QCoh(\LocSys) 
\overset{\Gamma(\LocSys,-)}\longrightarrow \Vect,
\end{multline*} 
while the composition
\begin{multline*} 
\Vect \to D(\Bun_G)\underset{\QCoh(\LocSys)}\otimes D(\Bun_G) \simeq \\
\simeq \left(D(\Bun_G)\otimes D(\Bun_G)\right)\underset{\QCoh(\LocSys\times \LocSys)}\otimes \QCoh(\LocSys)
\to D(\Bun_G)\otimes D(\Bun_G)
\end{multline*}
is the unit of the absolute self-duality of $D(\Bun_G)$, i.e.,
$$\Vect \overset{\Qlb\mapsto \omega_{\Bun_G}}\longrightarrow D(\Bun_G)\overset{(\Delta_{\Bun_G})_\bullet}\longrightarrow D(\Bun_G\times \Bun_G)\simeq 
D(\Bun_G)\otimes D(\Bun_G).$$

\medskip

Hence, the map in \eqref{e:long composition} identifies with
\begin{multline}    \label{e:long composition rewrite}
\Vect \overset{\Qlb\mapsto \omega_{\Bun_G}}\longrightarrow D(\Bun_G)\overset{(\Delta_{\Bun_G})_\bullet}\longrightarrow D(\Bun_G\times \Bun_G) \simeq \\
\simeq D(\Bun_G)\otimes D(\Bun_G) \overset{\Phi}\longrightarrow 
\QCoh(\LocSys)\overset{\otimes \CE_{V,x}}\longrightarrow \QCoh(\LocSys) 
\overset{\Gamma(\LocSys,-)}\longrightarrow \Vect.
\end{multline} 

\sssec{}

Using \eqref{e:interp counit} and using the fact that the functor $\CE_{V,x}\star -$ on $D(\Bun_G)$ is the Hecke functor $H_{V,x}$,  
we rewrite the map in \eqref{e:long composition rewrite} as 
\begin{multline*} 
\Vect \overset{\Qlb\mapsto \omega_{\Bun_G}}\longrightarrow D(\Bun_G)\overset{(\Delta_{\Bun_G})_\bullet}\longrightarrow D(\Bun_G\times \Bun_G)  \simeq \\
\simeq D(\Bun_G)\otimes D(\Bun_G) \overset{\on{Id}\otimes H_{V,x}} \longrightarrow D(\Bun_G)\otimes D(\Bun_G) 
\overset{\Frob^!_{\Bun_G}\otimes \on{Id}} \longrightarrow D(\Bun_G)\otimes D(\Bun_G) \to \\
\overset{(\Delta_{\Bun_G})^!} \longrightarrow 
D(\Bun_G) \overset{\on{H}^\cdot(\Bun_G,-)}\longrightarrow \Vect.
\end{multline*} 

However, by base change, the latter is isomorphic to 
$$\on{H}^\cdot(\Bun_G\underset{\Frob_{\Bun_G},\Bun_G\times \Bun_G}\times \CH_x,\CS_V),$$
as required.

\section{Higher representations vs. usual representations}

\ssec{Representations of groups on categories}  \label{ss:higher rep}

\sssec{}

Let $H$ be a group (ind-scheme) over $\Fqb$. We consider the 2-category 
$$H\mmod:=D(H)\mmod$$ as an object of $\tDGCat$.

\sssec{}  \label{sss:tw Frob}

Assume that $H$ is defined over $\Fq$. Then the Frobenius endomorphism of $H$ defines an endo-functor
$\Frob^!_H$ of $H\mmod$. 

\medskip

We claim that
\begin{equation} \label{e:Tr on rep}
\Tr(\Frob^!_H,H\mmod)=D(H)_{\on{Ad}_{\Frob}(H)},
\end{equation}
where the subscript $\on{Ad}_{\Frob}(H)$ means conviariants with respect to the Frobenius-twisted conjugation.

\medskip

Indeed, 
$$\Tr(\Frob^!_{H},H\mmod)\simeq D(H)\underset{\on{Graph}_{\Frob_{H}},D(H)\otimes D(H),\Delta_{H}}\otimes D(H),$$
which can be rewritten as
$$D(H)\underset{D(H)}\otimes \Vect,$$
where $D(H)$ acts on itself by Frobenius-twisted conjugation:
$$h(h')=\Frob_{H}(h)\cdot h'\cdot h^{-1}.$$

\ssec{Representations of Chevalley groups}

\sssec{}

There are two main cases that we want to consider: one is when $H$ is an algebraic group $G$, considered
in this section, and the other is when $H$ is loop group $G(\CK)$, where $\CK$ is a local field and $G$ is 
reductive (considered in the next section). 

\medskip

For $H=G$, we will pretend that $G\mmod$ behaves as in the case of D-modules (and hopefully, Drinfeld's
formalism will justify that). 

\medskip

Namely, that for $\CC\in G\mmod$, the *-averaging functor defines an equivalence
\begin{equation} \label{e:coinv}
\CC_G\to \CC^G.
\end{equation} 

\sssec{}

Thus, in particular, we identify 
$$D(G)_{\on{Ad}_{\Frob}(G)}\simeq D(G)^{\on{Ad}_{\Frob}(G)}.$$

\medskip

Assume now that $G$ is connected. Then by Lang's theorem, the $\on{Ad}_{\Frob}(G)$-action of $G$ 
on itself is transitive and the stabilizer of $1\in G$ is the finite group $G(\Fq)$. 
Hence, we obtain
$$D(G)^{\on{Ad}_{\Frob}(G)}\simeq \Rep(G(\Fq)).$$

Thus, we obtain a canonical equivalence of categories:

\begin{equation} \label{e:rep arise}
\Tr(\Frob^!_G,G\mmod)=\Rep(G(\Fq)).
\end{equation}

I.e., the category of representations of the Chevalley group $G(\Fq)$ arises from the 2-category $G\mmod$ as the trace of
Frobenuis. 

\sssec{}

Let $\CC$ be a dualizable object of $G\mmod$. The equation \eqref{e:coinv} implies that this condition is equivalent to being
dualizable as a plain DG-category. Moreover, the dual of $\CC$, viewed as a category, equipped with an action of $G$, identifies
canonically with $\CC^\vee$, equipped with the natural $G$-action. 

\medskip

Let $\CC$ be equipped with a functor 
\begin{equation} \label{e:Frob on cat}
\Frob_\CC:\CC\to \CC
\end{equation}
compatible with the monoidal endomorphism $\Frob^!_{G}$ of $D(G)$. Then, following \secref{sss:higher trace}
we attach to $(\CC,\Frob_\CC)$ an object 
$$\Tr_{\Frob^!_{G},D(G)}(\Frob_\CC,\CC)\in \Rep(G(\Fq)).$$

\medskip

By \secref{sss:comp}, the vector space underlying $\Tr_{\Frob^!_{G},D(G)}(\Frob_\CC,\CC)$ is simply
$\Tr(\Frob_\CC,\CC)$.

\ssec{An example: Deligne-Lusztig representations}

\sssec{}  \label{sss:DL}

Let $G$ be reductive, and take $\CC=D(G/B)$. We let the datum of \eqref{e:Frob on cat} be
the composition of the usual Frobenius on $G/B$ and the pull-push along the correspondence
\begin{equation} \label{e:pull-push}
G/B  \leftarrow (G/B\times G/B)_w\to G/B,
\end{equation} 
where $(G/B\times G/B)_w\subset (G/B\times G/B)$ be the subvariety of pairs of Borels in relative position $w\in W$.

\medskip

Then the corresponding object 
$$\Tr_{\Frob^!_{G},D(G)}(\Frob_\CC,\CC)\in \Rep(G(\Fq))$$
is by definition the Deligne-Lusztig representation corresponding to $w$.

\sssec{}

Note that the pull-push functor in \eqref{e:pull-push} can also be interpreted as the convolution on the right
with the $w$-costandard object on $B\backslash G/B$. In particular, it is an equivalence, and admits a right
adjoint, which is also given by a kernel (convolution on the right with the $w^{-1}$-standard). 

\ssec{Traces}

\sssec{}

Consider the trace of the \emph{identity} endomorphism of $G\mmod$. In a way similar to \secref{sss:tw Frob},
the resulting category identifies with 
$$D(G)_{\on{Ad}(G)},$$
and using \eqref{e:coinv}, further with 
$$D(G)^{\on{Ad}_(G)}.$$

\medskip

To a dualizable $\CC\in G\mmod$, we can thus attach an object 
$$\chi(\CC)\in D(G)^{\on{Ad}(G)}.$$

\medskip

By construction, the !-fiber of $\chi(\CC)$ at $g\in G$ equals the trace of the endo-functor of $\CC$, given by $g$.

\sssec{}

For example, for $\CC=D(G/B)$,
$$\chi(\CC)=\on{Spr},$$
where $\on{Spr}$ denotes the Springer sheaf. 

\sssec{}  

Suppose again that $\CC$ is endowed with a datum of \eqref{e:Frob on cat}. Assume, moreover, that the functor
$\Frob_\CC$ admits a right adjoint. Then the construction from \secref{sss:trace of 1-morph} gives rise
to a map
\begin{equation} \label{e:Weil structure}
\chi(\CC)\to \Frob^!(\chi(\CC)).
\end{equation} 

\sssec{Example}

One can show that for $(\CC,\Frob_\CC)$ as in \secref{sss:DL}, the resulting map \eqref{e:Weil structure}
is given by the composing the natural Weil structure on $\on{Spr}$ and the automorphism of $\on{Spr}$,
given by the element $w\in W$. This is so, because this is obviously so over the locus of $G$, consisting
of regular semi-simple elements. 

\sssec{}

By imposing certain finiteness conditions on $\CC$ (in the spirit of \emph{complete dualizability}) one can 
ensure that the objects 
$$\chi(\CC)\in D(G)^{\on{Ad}(G)} \text{ and } \Tr_{\Frob^!_{G},D(G)}(\Frob_\CC,\CC)\in \Rep(G(\Fq))$$
are both compact.

\medskip

In this case one can show that 
$$\Tr(\Frob_G^!,\chi(\CC)),$$
which is an Ad-invariant function on $G(\Fq)$ equals the character of the representation $\Tr_{\Frob^!_{G},D(G)}(\Frob_\CC,\CC)$.

\sssec{}

In particular, this shows that the character of the Deligne-Lusztig representation corresponding to $w$ is equal to the function
obtained from $\on{Spr}$, by twisting its Weil structure by $w$. 

\section{The local story}

\ssec{Categories acted on by the loop group}

\sssec{}

We now consider the situation of \secref{ss:higher rep} for $H=G(\CK)$.

\medskip

In this case (at least in the case of D-modules), we still have an equivalence
\begin{equation} \label{e:loop coinv}
\CC_{G(\CK)}\to \CC^{G(\CK)}.
\end{equation}

But this equivalence depends on the choice of a point of a parahoric subgroup $P$. We make this choice to be $P=G(\CO)$.

\sssec{}

Namely, the equivalence
\begin{equation} \label{e:arc coinv}
\CC_{P}\to \CC^{P}
\end{equation}
is the same as in \eqref{e:coinv}. The category $\CC_{P}$ is acted by the Hecke category $\CH:=D(P\backslash G(\CK)/P)$.
We interpret
$$\CC^{G(\CK)}\simeq \on{Funct}_{\CH}(D(\on{pt}/P),\CC^{P}),$$
and using \eqref{e:arc coinv} also 
$$\CC_{G(\CK)}\simeq \CC^{P}\underset{\CH}\otimes D(\on{pt}/P).$$

Now, the \emph{ind-properness} of $G(\CK)/P$ implies that for any $\CC'\in \CH\mmod$, we have a canonical equivalence
$$\on{Funct}_{\CH}(D(\on{pt}/P),\CC')\simeq \CC'\underset{\CH}\otimes D(\on{pt}/P),$$
characterized by the fact that the natural projection
$$\CC'\simeq \CC' \underset{D(\on{pt}/P)}\otimes D(\on{pt}/P)\to \CC'\underset{\CH}\otimes D(\on{pt}/P)$$
becomes identified with the \emph{left} adjoint of the forgetful functor
$$\on{Funct}_{\CH}(D(\on{pt}/P),\CC')\to \on{Funct}_{D(\on{pt}/P)}(D(\on{pt}/P),\CC')\simeq \CC'.$$

\sssec{}

If we choose a different parabolic, say $P'\subset P$, the two identifications in \eqref{e:loop coinv} will differ by $\det(H^\cdot_c(P/P',k))$. 
I.e., the choice in \eqref{e:loop coinv} is really that of a trivialization of the dimension torsor on $G(\CK)$. 

\ssec{The local category}

\sssec{}

We now assume that the ground field is $\Fqb$ and that $G$ is defined over $\Fq$. Combing \eqref{e:loop coinv} with
\eqref{e:Tr on rep}, we obtain an identification 
\begin{equation} \label{e:loc category}
\Tr(\Frob^!_{G(\CK)},G(\CK)\mmod)\simeq D(G(\CK))^{\on{Ad}_{\Frob}(G(\CK))}.
\end{equation} 

The category in the right-hand side in \eqref{e:loc category} is \emph{the} category on the geometric side of Langlands
that we propose to study in a joint project with Alain Genestier and Vincent Lafforgue. 

\medskip

This category is related to one appearing in Fargues' conjecture \cite{Far}.  

\sssec{}

Taking the !-fiber at $1\in G(\CK)$ we obtain a functor
$$\sR:D(G(\CK)^{\on{Ad}_{\Frob}(G(\CK))}\to \Rep(G(\CK)(\Fq)).$$

However, since Lang's theorem fails for $G(\CK)$, the above functor $\sR$ is no longer an equivalence.

\sssec{}

Let $\CC$ be a dualizable object of $G(\CK)\mmod$. As in the case of a finite-dimensional $G$, the equivalence
\eqref{e:loop coinv} implies that the above dualizability condition is equivalent to $\CC$ being dualizable as a
plain object of $\DGCat$. 

\medskip

Moreover, the dual of $\CC$ as a category acted on by $G(\CK)$ identifies with $\CC^\vee$,
equipped with the natural $G(\CK)$-action (but this identification depends on the choice of a trivialization of the dimension torsor on $G(\CK)$).  

\sssec{}  \label{sss:loc Categ tr}

Assume that $\CC$ is equipped with a functor $$\Frob_\CC:\CC\to \CC$$ compatible 
with the monoidal endomorphism $\Frob^!_{G(\CK)}$ of $D(G(\CK))$.

\medskip

We obtain that to this datum one can canonically attach an object
$$\Tr_{\Frob^!_{G(\CK)},D(G(\CK))}(\Frob_\CC,\CC)\in D(G(\CK))^{\on{Ad}_{\Frob}(G(\CK))}.$$

\medskip

The !-fiber of $\Tr_{\Frob^!_{G(\CK)},D(G(\CK))}(\Frob_\CC,\CC)$ at $1\in G(\CK)$, viewed as a mere vector space identifies
with 
$$\Tr(\Frob_\CC,\CC).$$
This identification also depends on the choice of a trivialization of the dimension torsor on $G(\CK)$. 

\medskip

The group $G(\CK)(\Fq)$ acts on $\Tr(\Frob_\CC,\CC)$ by transport of structure. This action coincides with one
obtained by considering the object
$$\sR\left(\Tr_{\Frob^!_{G(\CK)},D(G(\CK))}(\Frob_\CC,\CC)\right).$$

\ssec{Action of local systems}

\sssec{}

Let $\LocSys_{\on{loc}}$ be the stack of local systems on the punctured disc. The local version of Hecke action says that
the algebra object $\QCoh(\LocSys_{\on{loc}})\mmod$ in $\tDGCat$ acts on $G(\CK)\mmod$.

\medskip

This structure can be also formulated as saying that 
the symmetric monoidal category $\QCoh(\LocSys_{\on{loc}})$ acts canonically on every object $\CC\in G(\CK)\mmod$. 

\sssec{}

By \secref{sss:trace of 1-morph} we obtain an action of the monoidal category
$$\Tr(\Frob^*_{\LocSys_{\on{loc}}},\QCoh(\LocSys_{\on{loc}})\mmod)\simeq \QCoh(\LocSys_{\on{loc}}^{\on{arthm}})$$
on $$\Tr(\Frob^!_{G(\CK)},G(\CK)\mmod)=D(G(\CK))^{\on{Ad}_{\Frob}(G(\CK))}.$$

\medskip

This is the action that we aim to construct in a joint project with Alain and Vincent.

\sssec{}

Let $(G(\CK)\mmod)_{\on{temp}}$ be the tempered part of $G(\CK)\mmod$. We expect that the 2-category $(G(\CK)\mmod)_{\on{temp}}$, when viewed
as a module over $\QCoh(\LocSys_{\on{loc}})$, is free on one generator.

\medskip

Define
$$D(G(\CK))^{\on{Ad}_{\Frob}(G(\CK))}_{\on{temp}}: =\Tr(\Frob^!_{G(\CK)},G(\CK)\mmod_{\on{temp}}).$$

We obtain that the category $D(G(\CK))^{\on{Ad}_{\Frob}(G(\CK))}_{\on{temp}}$ is a free module on one generator over 
$\QCoh(\LocSys_{\on{loc}}^{\on{arthm}})$.

\ssec{Imposing $G(\CO)$-invariance}

\sssec{}

It is an expectation of the local \emph{geometric} Langlands that the 2-category
$$G(\CK)\mmod \underset{\QCoh(\LocSys_{\on{loc}})\mmod}\otimes \QCoh(\LocSys_{\on{loc-unr}})\mmod$$
is equivalent to $\CH\mmod$, where
$$\CH=D(G(\CO)\backslash G(\CK)/G(\CO))$$
is the spherical Hecke category.

\medskip

In the above formula $\LocSys_{\on{loc-unr}}$ is the stack of unramified geometric local systems on the formal disc around $x$
(i.e., the classifying stack $\on{pt}/\cG$).

\begin{rem}
Similarly, we expect that 
$$G(\CK)\mmod \underset{\QCoh(\LocSys_{\on{loc}})\mmod}\otimes \QCoh(\LocSys_{\on{loc-tame}})\mmod$$
is equivalent to the 2-category of module categories over the Iwahori-Hecke category.
\end{rem}

\sssec{}

Taking the traces of the Frobenius (here I'm assuming that the operation of taking the trace of an endo-functor
is compatible with geometric realizations under certain hypothesis), we obtain that the category
$$\Tr(\Frob^!_{G(\CK)},G(\CK)\mmod)\underset{\QCoh(\LocSys^{\on{arthm}}_{\on{loc}})}\otimes \QCoh(\LocSys^{\on{arthm}}_{\on{loc-unr}})$$
is canonically equivalent to the Hochschild homology of $\CH$ regarded as a bi-module category over itself, with the left action being
the monoidal operation and the right action is twisted by $\Frob^!$.

\medskip

One shows that the latter Hochshild homology identifies with $\QCoh(\cG/\on{Ad}(\cG))$. Indeed, this would be obviously so if instead
of $\CH$ we took the naive Hecke category, i.e., $\Rep(\cG)$. Now, the quasi-coherent interpretation of $\CH$ shows that the derived
stuff disappears under the Frobenius (note, however, that a parallel fact fails completely if instead of the Frobenius we consider the
identity functor, i.e., when we consider the plain Hochschild homology of $\CH$). 

\medskip

Combining with \eqref{e:loc category}, we obtain an equivalence
\begin{equation} \label{e:local unramified}
D(G(\CK))^{\on{Ad}_{\Frob}(G(\CK))}\underset{\QCoh(\LocSys^{\on{arthm}}_{\on{loc}})}\otimes \QCoh(\LocSys^{\on{arthm}}_{\on{loc-unr}})
\simeq \QCoh(\cG/\on{Ad}(\cG)).
\end{equation}

\sssec{}  \label{sss:G(O)-inv}

Let $(\CC,\Frob_\CC)$ be as in \secref{sss:loc Categ tr}. On the one hand, consider the image of 
$\Tr_{\Frob^!_{G(\CK)},D(G(\CK))}(\Frob_\CC,\CC)$ 
under the functor
$$D(G(\CK))^{\on{Ad}_{\Frob}(G(\CK))}\to 
D(G(\CK))^{\on{Ad}_{\Frob}(G(\CK))}\underset{\QCoh(\LocSys^{\on{arthm}}_{\on{loc}})}\otimes \QCoh(\LocSys^{\on{arthm}}_{\on{loc-unr}}).$$

On the other hand, consider $\CC^{G(\CO)}$ as an object of $\CH\mmod$, and consider
$$\Tr_{\Frob^!,\CH}(\Frob_\CC,\CC^{G(\CO)})\in \QCoh(\cG/\on{Ad}(\cG)).$$

It follows from the definitions, that the above two objects are identified under \eqref{e:local unramified}.

\ssec{Local vs global compatibilty}

\sssec{}

Let $X$ be again a complete curve, and $x\in X$ a rational point. Consider $\Bun_G^{\on{level}_x}$ as a stack equipped with an action of
$G(\CK)$, where $\CK$ is the local field at $x$. So $D(\Bun_G^{\on{level}_x})$ is naturally an object of $G(\CK)\mmod$, equipped with
a compatible Frobenius. 

\medskip

Hence, by \secref{sss:loc Categ tr}, we have a canonically defined object
\begin{equation} \label{e:Drinf loc coarse}
\Tr_{\Frob^!_{G(\CK)},D(G(\CK))}(\Frob^!_{\Bun_G^{\on{level}_x}},D(\Bun_G^{\on{level}_x}))\in D(G(\CK))^{\on{Ad}_{\Frob}(G(\CK))}.
\end{equation} 

\medskip

By construction, we have
$$\sR\left(\Tr_{\Frob^!_{G(\CK)},D(G(\CK))}(\Frob^!_{\Bun_G^{\on{level}_x}},D(\Bun_G^{\on{level}_x}))\right)\simeq \on{Funct}(\Bun_G^{\on{level}_x}(\Fq),\Qlb),$$
as representations of $G(\CK)(\Fq)$. 

\sssec{}

Let $\LocSys_{\on{glob}}$ denote the stack of local systems on $X-x$. We have the geometric Hecke action of the monoidal category
$\QCoh(\LocSys_{\on{glob}})$ on $D(\Bun_G^{\on{level}_x})$. 

\medskip

As in \secref{ss:act loc sys}, this action defines an action of the commutative algebra 
$$\Gamma(\LocSys^{\on{arthm}}_{\on{glob}},\CO_{\LocSys^{\on{arthm}}_{\on{glob}}})$$
on the object \eqref{e:Drinf loc coarse}. 

\sssec{}

In particular, by functoriality, we obtain an action of $\Gamma(\LocSys^{\on{arthm}}_{\on{glob}},\CO_{\LocSys^{\on{arthm}}_{\on{glob}}})$
on $\on{Funct}(\Bun_G^{\on{level}_x}(\Fq),\Qlb)$.

\medskip

However, the above action of $\Gamma(\LocSys^{\on{arthm}}_{\on{glob}},\CO_{\LocSys^{\on{arthm}}_{\on{glob}}})$ on $\on{Funct}(\Bun_G^{\on{level}_x}(\Fq),\Qlb)$
produces nothing essentially new as compared to the situation of \secref{ss:act loc sys}: this is the action obtained from the action of
$\QCoh(\LocSys_{\on{glob}})$ on $D(\Bun_G^{\on{level}_x})$ by taking traces of the Frobenius endo-morphisms.

\sssec{}

Here is, however, an object that did not appear previously: the procedure from \secref{ss:Drinf} allows to refine the object \eqref{e:Drinf loc coarse}
to an object 

$$\Tr_{\Frob^!_{G(\CK)},D(G(\CK))}(\Frob^!_{\Bun_G^{\on{level}_x}},D(\Bun_G^{\on{level}_x}))\in D(G(\CK))^{\on{Ad}_{\Frob}(G(\CK))}$$
to an object 
$$\on{Drinf}_x\in D(G(\CK))^{\on{Ad}_{\Frob}(G(\CK))}\underset{\QCoh(\LocSys^{\on{arthm}}_{\on{loc}})}\otimes \QCoh(\LocSys^{\on{arthm}}_{\on{glob}}).$$

\sssec{}

Here is how the above object $\on{Drinf}_x$ is related to 
$$\on{Drinf}\in \QCoh(\LocSys^{\on{arthm}}_{\on{glob-unr}})$$
from \secref{ss:Drinf}.

\medskip

Consider the image of $\on{Drinf}_x$ under the functor
\begin{multline*}
D(G(\CK))^{\on{Ad}_{\Frob}(G(\CK))}\underset{\QCoh(\LocSys^{\on{arthm}}_{\on{loc}})}\otimes \QCoh(\LocSys^{\on{arthm}}_{\on{glob}})\to \\
\to D(G(\CK))^{\on{Ad}_{\Frob}(G(\CK))}\underset{\QCoh(\LocSys^{\on{arthm}}_{\on{loc}})}\otimes \QCoh(\LocSys^{\on{arthm}}_{\on{glob}})
\underset{\QCoh(\LocSys^{\on{arthm}}_{\on{glob}})}\otimes \QCoh(\LocSys^{\on{arthm}}_{\on{glob-unr}}).
\end{multline*}

Note, however, that the latter category identifies with
$$D(G(\CK))^{\on{Ad}_{\Frob}(G(\CK))}\underset{\QCoh(\LocSys^{\on{arthm}}_{\on{loc}})}\otimes \QCoh(\LocSys^{\on{arthm}}_{\on{loc-unr}})
\underset{\QCoh(\LocSys^{\on{arthm}}_{\on{loc-unr}})}\otimes \QCoh(\LocSys^{\on{arthm}}_{\on{glob-unr}}),$$
and further by \eqref{e:local unramified} with
\begin{multline*}
\QCoh(\cG/\on{Ad}(\cG)) \underset{\QCoh(\LocSys^{\on{arthm}}_{\on{loc-unr}})}\otimes 
\QCoh( \LocSys^{\on{arthm}}_{\on{glob-unr}})\simeq \\
\simeq \QCoh(\cG/\on{Ad}(\cG)) \underset{\QCoh(\cG/\on{Ad}(\cG))}\otimes 
\QCoh( \LocSys^{\on{arthm}}_{\on{glob-unr}})\simeq \QCoh(\LocSys^{\on{arthm}}_{\on{glob-unr}}).
\end{multline*}

\medskip

Thus, we obtain a functor
$$D(G(\CK))^{\on{Ad}_{\Frob}(G(\CK))}\underset{\QCoh(\LocSys^{\on{arthm}}_{\on{loc}})}\otimes \QCoh(\LocSys^{\on{arthm}}_{\on{glob}})\to
\QCoh(\LocSys^{\on{arthm}}_{\on{glob-unr}}).$$

Now, it follows from \secref{sss:G(O)-inv} that the image of $\on{Drinf}_x$ under the above functor identifies with $\on{Drinf}$.


\begin{thebibliography}{99}

\bibitem[Ga]{Ga} D.~Gaitsgory, {\it Outline of the proof of the geometric Langlands conjecture for $GL(2)$}, arXiv:1203.6343. 

\bibitem[Far]{Far} L.~Fargues, {\it Geometrization of the local Langlands correspondence: an overview}, \hfill \newline
https://webusers.imj-prg.fr/$~$laurent.fargues/Prepublications.html

\bibitem[Laf]{Laf} V.~Lafforgue, {\it Chtoucas pour les groupes r\'eductifs et param\'etrisation de Langlands globale}, arXiv:1209.5352.


\end{thebibliography}
\end{document}